\documentclass[11pt]{amsart}
\usepackage{amssymb}
\usepackage{amsmath}
\usepackage{amssymb}
\usepackage{amsfonts}
\usepackage{tipa}
\usepackage{amssymb}
\usepackage{mathrsfs}
\usepackage{amsmath}
\usepackage{txfonts}

\theoremstyle{plain}
\newtheorem{Thm}{Theorem}

\newtheorem{Prop}[Thm]{Proposition}

\begin{document}
\large

\title[Li-Yau-Hamilton type estimates for heat equations]
{Li-Yau-Hamilton type estimates for heat equations on manifolds}

\author{Li Ma}

\address{Department of mathematics \\
Henan Normal university \\
Xinxiang, 453007 \\
China}

\email{nuslma@gmail.com}

\thanks{The research is partially supported by the National Natural Science
Foundation of China 1127111 and SRFDP 20090002110019}

\begin{abstract}
In this paper, we study the gradient estimates of Li-Yau-Hamilton
type for positive solutions to both drifting heat equation and the
simple nonlinear heat equation problem
$$
u_t-\Delta u=au\log u, \ \  u>0
$$
on  the compact Riemannian manifold $(M,g)$ of dimension $n$ and
with non-negative (Bakry-Emery)-Ricci curvature. Here $a\leq 0$ is
a constant. The latter heat equation is a basic evolution equation
which is the negative gradient heat flow to the functional of
Log-Sobolev inequality on the Riemannian manifold. We derive various versions of gradient estimates which generalize
Hamilton's gradient estimate. An
question concerning the Hamilton type gradient estimate for the simple nonlinear heat equation is
addressed.


\end{abstract}

 \maketitle

 \newpage

\section{Introduction}

In this paper, we generalize the previous gradient estimate of Li-Yau-Hamilton type
\cite{H} for heat equations, which have been studied recently by
many researchers and can be seen for example, \cite{P02}, \cite{Chow}, \cite{Ni},\cite{BCP},
\cite{HM}, \cite{LX}, \cite{L},\cite{M2},\cite{KZ},
\cite{CZ},\cite{CH},\cite{Sun},\cite{CC},\cite{qian}, etc. The
equations under consideration have their deep background from the
fundamental gap of the Schrodinger equation and the Ricci flow (see
\cite{Hs}, \cite{ML}, \cite{Yst}, etc). Interesting gradient
estimates of Hamilton type for heat equation associated to Ricci
flow have been obtained in \cite{CH}, \cite{BCP}, and \cite{Z}. We
derive the Hamilton type gradient estimate for the drifting heat
equation and the simple nonlinear heat equation from the view-point
of the Bernstein type estimates. This is a new observation and our
argument is shorter than previous ones. This paper is a continuation of our work \cite{MA}.

We first derive the Li-Yau-Hamilton type gradient estimate for the
drifting heat equation
\begin{equation}\label{ht}
u_t-\Delta u=-\nabla\phi\cdot\nabla u, \ \  u>0,  t\in (0,T]
\end{equation}
on the compact Riemannian manifold $(M,g)$ of dimension $n$. Here $T>0$,
$\phi$ is a smooth function on $M$, and $\nabla \phi$ is the
gradient of $\phi$ in the metric $g$. We shall denote $D^2\phi$ the
hessian matrix of $\phi$.

We have the following gradient estimate for (\ref{ht}) with $\phi=0$, which is of the type of Hamilton's gradient estimate.
\begin{Thm}\label{main}
Assume that the compact Riemannian manifold $(M,g)$ has the
non-negative Ricci curvature. Let $u>0$ be a positive smooth
solution to (\ref{ht}) with $\phi=0$. Assume that $\sup_Mu=1$. Then we have for any $p\in (0,1)$,
$$
tu^{2p-2}|\nabla u|^2\leq \frac{1-u^p}{p(1-p)}, \ \ in \ \ M\times [0,T].
$$
The same estimate is true for (\ref{ht}) with $\phi=0$ on the compact Riemannian manifold with Neumann boundary condition or on complete Riemannian
manifolds when the maximum principle can be applied.
\end{Thm}
This result can be generalized to the case for the drifting heat equation (\ref{ht}) with non-trivial $\phi$ and with nongative Bakry-Emory-Ricci curvature. For simply written purpose, we just present the case when $p=1$ below.

\begin{Thm}\label{thm1}
Assume that the compact Riemannian manifold $(M,g)$ has the
non-negative Bakry-Emery-Ricci curvature in the sense that
$$Rc+D^2\phi\geq -K
$$ on $M$ for some constant $K\geq 0$. Let $u>0$ be a positive smooth
solution to (\ref{ht}). Assume that $\sup_Mu=1$. Let $f=-\log u$.
Then we have, for all $t>0$,
$$
t|\nabla f|^2\leq (2Kt+1)f.
$$
The same estimate is true for (\ref{ht}) on complete Riemannian
manifolds when the maximum principle can be applied.
\end{Thm}

We also have the following result for (\ref{ht}) on a compact Riemannian manifold
with smooth boundary.

\begin{Thm}\label{thm2}
Assume that the compact Riemannian manifold $(M,g)$ with convex
boundary has the curvature condition about the Bakry-Emery-Ricci
tensor that $$Rc+D^2\phi\geq -K$$ on $M$ for some constant $K\geq
0$. Let $u>0$ be a positive smooth solution to (\ref{ht}) with
Neumann boundary condition $u_\nu=0$, where $\nu$ is the outward
unit normal to the boundary. Assume that $\sup_Mu=1$. Let $v=\sqrt{u}$ and $f=-\log
u$. Then we have, for all $t>0$,
\begin{equation}\label{est-1}
t|\nabla v|^2\leq (2Kt+1)fu.
\end{equation}
In case when $Rc\geq 0$ on $M$, we also have
we have the gradient estimate
$$
t|\nabla u|^2\leq \frac{u^{2-p}}{p(1-p)}, \ \ in \ \ M\times [0,T]
$$
for any $p\in [1/2, 1)$.
\end{Thm}

Assume that $K=0$ and $u$ is any bounded smooth solution to
(\ref{ht}). Assume that $A=\sup_Mu>0$. Let $v=(A-u)/A$. Then $v$
is a positive solution to (\ref{ht}) and the above gradient
estimate is
\begin{equation}\label{eq1-1}
 t|\nabla u|^2\leq (A-u)^2\log \frac{A}{A-u},
\end{equation}
which is the usual form of Hamilton type gradient estimate.

The drifting heat equation is closely related to the fundamental
gap of the Schrodinger operator on convex domains (so $K=0$).
Namely, Let $\lambda =\lambda_2-\lambda_1$ be the fundamental gap
of the Laplacian operator $-\Delta$ and let $f_j$ be the
eigenfunctions corresponding to $\lambda_j$, $j=1,2$. Let
$u:=u(x):=f_2/f_1$. Then we have (\cite{Yst})
$$
\Delta u=-\lambda u-2(\nabla u\cdot \nabla \log f_1).
$$
Set
$$
\phi=-2\log f_1,
$$
which is convex by the well-known result of Brascamp-Lieb
\cite{BL}. Let
$$
v(x,t)=exp(-\lambda t)u(x).
$$
Then $v$ satisfies (\ref{ht}) with the Neumann boundary condition.

The other interesting problem to us is to derive the Li-Yau-Hamilton
type gradient estimate for the following nonlinear heat equation
problem
\begin{equation}\label{n-ht1}
u_t-\Delta u=au\log u, \ \  u>0
\end{equation}
on  the compact Riemannian manifold $(M,g)$ of dimension $n$. Here
$a\in \mathbf{R}$ is some constant. This heat equation can be
considered as the negative gradient heat flow to $W$-functional
\cite{P02}, which is closely related to the Log-Sobolev
inequalities on the Riemannian manifold. In \cite{M}, we propose
the study of the local gradient estimates for solutions to
(\ref{n-ht1}) based on its relation with Ricci solitons. Soon
after, Y.Yang has studied this problem in \cite{Y} and his result is
Li-Yau type gradient estimate \cite{LY}.

We have the following result for (\ref{n-ht1}).
\begin{Thm}\label{thm3}
Assume that the compact Riemannian manifold $(M,g)$ has the
non-negative Ricci curvature condition, i.e., $Rc\geq 0$. Let
$u>0$ be a positive smooth solution to (\ref{n-ht1}). Assume that
$\sup_Mu<1$ at the initial time and $a\leq 0$. Let $f=-\log {u}$.
Then we have, for all $t>0$, $\sup_Mu<1$ and
$$
t|\nabla f|^2\leq f.
$$
The same estimate is true for (\ref{ht}) on complete Riemannian
manifolds when the maximum principle can be applied.
\end{Thm}

Similar to Theorem \ref{thm2}, we also have the following result
for (\ref{n-ht1}) on the manifold with smooth boundary.

\begin{Thm}\label{thm4}
Assume that the compact Riemannian manifold $(M,g)$ with convex
boundary has the non-negative Ricci curvature condition. Let $u>0$
be a positive smooth solution to (\ref{ht}) with Neumann boundary
condition $u_\nu=0$, where $\nu$ is the outward unit normal to the
boundary. Assume that $\sup_Mu<1$ at the initial time and $a\leq
0$. Let $f=-\log u$. Then we have, for all $t>0$, $\sup_Mu<1$ and
\begin{equation}\label{est-2}
t|\nabla f|^2\leq f. \end{equation}
\end{Thm}

It is quite clear that our results to (\ref{n-ht1}) are not
satisfied by us because of the assumption $\sup_Mu<1$ at the initial
time. So we leave it open to derive the Hamilton type gradient
estimate for positive solutions to (\ref{n-ht1}). We remark that in
the recent work \cite{Hi}, the existence of positive solution to
stationary version of (\ref{n-ht1}) on a complete Riemannian
manifold is discussed.

The plan of our paper is below. In section \ref{sect2}, we give
the proofs of Theorems \ref{thm1} and \ref{thm2}. In section
\ref{sect3} we study (\ref{n-ht1}).

\section{Hamilton type estimate for the heat
equation}\label{sect1}
We assume that $Rc\geq 0$ on the compact Riemannian manifold $(M,g)$. Let $T>0$.

Assume that $u>0$ is a bounded positive solution to the heat equation, i.e., (\ref{ht}) with $\phi=0$ on $M\times [0,T]$. We may assume that $0<u\leq 1=\sup u$. Let
$f=-\log u$. Then
$$
f_j=-u_j/u, \ \ \Delta f=-\Delta u/u+|\nabla f|^2.
$$
Then we have
\begin{equation}\label{eq1-1}
(\partial_t-\Delta)f=-|\nabla f|^2.
\end{equation}
and
$$
(\partial_t-\Delta)(uf)=u|\nabla f|^2.
$$

Define $v=\sqrt{u}$. Then we have $u=v^2$,
$$
(\partial_t-\Delta)u=2v(\partial_t-\Delta)v-2|\nabla v|^2,
$$
which implies that
\begin{equation}\label{eq1-11}
(\partial_t-\Delta)v=v^{-1}|\nabla v|^2.
\end{equation}
One may compute that
$$
(\partial_t-\Delta) |\nabla v|^2= 2<\nabla
v, \nabla (\partial_t-\Delta)v>-2|D^2v|^2-2Rc(\nabla v, \nabla v).
$$
Then we have
\begin{equation}\label{eq1-2}
(\partial_t-\Delta) |\nabla v|^2\leq -2v^{-2}|\nabla v|^4+2v^{-1}<\nabla
v, \nabla |\nabla v|^2>-2|D^2v|^2\leq 0.
\end{equation}
Using (\ref{eq1-11}) and (\ref{eq1-2}) we get that
$$
(\partial_t-\Delta) (t4|\nabla v|^2-fu)\leq 4t(\partial_t-\Delta) (|\nabla v|^2)\leq 0.
$$
By the Maximum principle we know that
$$
t4|\nabla v|^2-fu\leq \sup_{t=0}(-fu)=0.
$$
This is Hamilton's gradient estimate for the heat equation.

We now generalize the above gradient estimate for the heat equation.
Let $0<p<1$ and let $v=u^p$. Then
$$
v_t=pu^{p-1}u_t, \ \ v_{ij}=pu^{p-1}u_{ij}+p(p-1)u^{p-2}u_iu_j,
$$
and
$$
(\partial_t-\Delta)v=-p(p-1)u^{p-2}|\nabla u|^2=\frac{1-p}{p}v^{-1}|\nabla v|^2.
$$
Then we have
\begin{equation}\label{eq1-22}
(\partial_t-\Delta) |\nabla v|^2\leq -2\frac{1-p}{p}v^{-2}|\nabla v|^4+2\frac{1-p}{p}v^{-1}<\nabla
v, \nabla |\nabla v|^2>-2|D^2v|^2
\end{equation}
Take $\frac{1}{2}\leq p<1$. Then we get $\frac{1-p}{p}\geq (\frac{1-p}{p})^2$ and
$$
(\partial_t-\Delta) |\nabla v|^2\leq 0.
$$
Hence, we have
$$
(\partial_t-\Delta) (t|\nabla v|^2-\frac{p}{1-p}v)\leq |\nabla v|^2-v^{-1}|\nabla v|^2\leq 0.
$$
By the maximum principle, we have
$$
t|\nabla v|^2-\frac{p}{1-p}v\leq 0,
$$
which is equivalent to the following inequality
$$
t|\nabla u|^2\leq \frac{u^{2-p}}{p(1-p)}, \ \ in \ \ M\times [0,T]
$$
Thus we have proved the following result.

\begin{Prop}\label{mali}
Let $(M,g)$ be a compact Riemannian manifold with non-negative Ricci curvature.
Let $u>0$ be a bounded positive solution to the heat equation such that $0<u\leq 1=\sup u$ on $M\times [0,T]$. Then we have the gradient estimate
$$
t|\nabla u|^2\leq \frac{u^{2-p}}{p(1-p)}, \ \ in \ \ M\times [0,T]
$$
for any $p\in [1/2,1)$.
\end{Prop}

We now prove the generalize version of Hamilton's gradient estimate, Theorem \ref{main}.

Let $H_0=t|\nabla v|^2+\frac{p}{1-p}v$. Then we have
$$
(\partial_t-\Delta) H_0\leq |\nabla v|^2+2t\frac{1-p}{p}v^{-1}<\nabla
v, \nabla |\nabla v|^2>+v^{-1}|\nabla v|^2:=Q_0.
$$
If $H_0$ takes its maximum on $M\times [0,T]$ at some point $(x_0,t_0)$ with $t_0>0$, then at $(x_0,t_0)$, we have
$(\partial_t-\Delta) H_0\geq 0$ and
$t\nabla |\nabla v|^2=-\frac{p}{1-p}\nabla v$. The latter implies that
$$
2t\frac{1-p}{p}v^{-1}<\nabla
v, \nabla |\nabla v|^2>=-2v^{-1}|\nabla v|^2
$$
and $Q_0=(1-v^{-1})|\nabla v|^2\leq 0$, which implies that $u=1$ at $(x_0,t_0)$. This is impossible by the strong maximum principle. Hence we have
$$
H_0\leq \sup_{t=0}H_0=\frac{p}{1-p},
$$
which is equivalent to
$$
tu^{2p-2}|\nabla u|^2\leq \frac{1-u^p}{p(1-p)}, \ \ in \ \ M\times [0,T].
$$
This completes the proof of Theorem \ref{main}.

Similar result can be stated for the compact Riemannian manifold with Neumann boundary condition.

\section{Li-Yau-Hamilton type estimate for drifting heat
equation}\label{sect2}

Assume that $u>0$ is a positive solution to (\ref{ht}). Let
$f=-\log u$. Then
$$
f_j=-u_j/u, \ \ \Delta f=-\Delta u/u+|\nabla f|^2.
$$
Then we have
\begin{equation}\label{eq2-1}
(\partial_t-\Delta)f+\nabla\phi\cdot \nabla f=-|\nabla f|^2,
\end{equation}
and
\begin{equation}\label{eq2-11}
(\partial_t-\Delta)(uf)+\nabla\phi\cdot \nabla (uf)=u|\nabla f|^2,
\end{equation}

Let $L=\partial_t-\Delta+\nabla\phi\cdot$ and $v=\sqrt{u}$.
Then as before, we have
$$
Lv=v=v^{-1}|\nabla v|^2.
$$
We compute $L|\nabla
v|^2$.

Note that
$$
(|\nabla v|^2)_t=2<\nabla v,\nabla v_t>.
$$
Recall the Bochner formula that
$$
\Delta|\nabla v|^2=2|D^2v|^2+<\nabla v,\nabla \Delta v>+2Rc(\nabla
v,\nabla v>.
$$

Then we have \begin{equation}\label{eq2-2}
 L|\nabla v|^2= 2<\nabla
v, \nabla Lv>-2|D^2v|^2-2(Rc+D^2\phi)(\nabla v, \nabla v).
\end{equation}

By the Ricci curvature bound assumption, we have
$$
L|\nabla v|^2\leq 2<\nabla
v, \nabla (v^{-1}|\nabla v|^2)>-2|D^2v|^2+2K |\nabla v|^2.
$$

Then we have
$$
L|\nabla v|^2\leq 2K |\nabla v|^2.
$$
Then we have
\begin{equation}\label{eq2-3}
L(4t|\nabla v|^2)\leq 4(1+2Kt) |\nabla v|^2.
\end{equation}

Using (\ref{eq2-11}), we get from (\ref{eq2-3}) that
\begin{equation}\label{eq2-4}
L(4t|\nabla v|^2-(2Kt+1)uf)\leq -2Kuf.
\end{equation}
We may re-write (\ref{eq2-4}) as
$$
L(4t|\nabla v|^2-(2Kt+1)uf)\leq \frac{2K}{2K+1}(4t|\nabla
v|^2-(2K+1)uf)-\frac{2K}{2K+1}(4t|\nabla v|^2).
$$
Then we have
$$
L(4t|\nabla v|^2-(2Kt+1)uf)\leq \frac{2K}{2K+1}(4t|\nabla
v|^2-(2K+1)uf).
$$
Note that
$$
\sup_{t=0}(4t|\nabla v|^2-(2Kt+1)uf)=0.
$$
Applying the Maximum principle we obtain that
$$
4t|\nabla v|^2-(2Kt+1)uf\leq 0.
$$
The proof of other part is similar to Proposition \ref{mali} and we omit the detail.
This completes the proof of Theorem \ref{thm1}.

We now prove Theorem \ref{thm2}. We need to treat the boundary
term. Note that $f_\nu=0$ on the boundary. Then on the boundary,
$$
[4t|\nabla v|^2-(2Kt+1)uf]_\nu=8tf_jf_{j\nu}=-8II(\nabla f,\nabla
f)\leq 0.
$$
Hence by the strong maximum principle we know that the maximum
point of $4t|\nabla f|^2-(2Kt+1)uf$ can not occur at the boundary
point and then we have
$$
4t|\nabla f|^2-(2Kt+1)f\leq 0.
$$
The proof of other part is similar to Proposition \ref{mali} and we omit the detail.
This completes the proof of Theorem \ref{thm2}.

\section{Li-Yau-Hamilton type estimate for the simple nonlinear heat
equation}\label{sect3}

As before, we let $f=-\log u$. Then we have
\begin{equation}\label{eq3-1}
(\partial_t-\Delta)f=af-|\nabla f|^2.
\end{equation}
Using $a\leq 0$ and the maximum principle we know that if
$\inf_Mf>0$ at the initial time, then it is always positive for
$t>0$.

Let $L:=\partial_t-\Delta$ in this section.
 Compute,
 \begin{equation}\label{eq3-2}
 L|\nabla f|^2= 2<\nabla
f, \nabla Lf>-2|D^2f|^2-2Rc(\nabla f, \nabla f).
\end{equation}
Then we have
$$
L|\nabla f|^2=2a|\nabla f|^2- 2<\nabla f, \nabla |\nabla
f|^2>-2|D^2f|^2-2Rc(\nabla f, \nabla f).
$$
Using the non-negative Ricci curvature assumption we have
$$
L|\nabla f|^2\leq 2a|\nabla f|^2- 2<\nabla f, \nabla |\nabla
f|^2>.
$$

Then
$$
L(t|\nabla f|^2)\leq (2at+1)|\nabla f|^2- 2<\nabla f, \nabla
(t|\nabla f|^2)>.
$$
Using (\ref{eq3-1}) we get that $$ L(t|\nabla f|^2-f)\leq
2at|\nabla f|^2-af- 2<\nabla f, \nabla (t|\nabla f|^2-f)>.
$$
Let $H=t|\nabla f|^2-f$. Then $f=t|\nabla f|^2-H$. Hence we have
$$ LH\leq at|\nabla
f|^2+aH- 2<\nabla f, \nabla H>.
$$

Using the assumption that $a\leq 0$, we obtain that
$$ LH\leq aH- 2<\nabla f, \nabla H>.
$$
Applying the maximum principle to $H$ we know that $H\leq 0$. That
is,
$$
t|\nabla f|^2-f\leq 0,
$$
which is the desired gradient estimate of Hamilton type. Then we
complete the proof of Theorem \ref{thm3}.

Using the same argument as in the proof of Theorem \ref{thm2}, we
can prove Theorem \ref{thm4}.

\emph{Acknowledgement}: The author would also like to thank IHES,
France for host and the K.C.Wong foundation for support in 2010.

\end{document}